\input amstex 
\documentstyle{amsppt}
\input bull-ppt
\keyedby{285/axm}

\topmatter
\cvol{26}
\cvolyear{1992}
\cmonth{April}
\cyear{1992}
\cvolno{2}
\cpgs{315-316}
\title A quasi-polynomial bound for the diameter\\of 
graphs of polyhedra \endtitle
\author Gil Kalai and Daniel J. Kleitman\endauthor
\shortauthor{Gil Kalai and D. J. Kleitman}
\shorttitle{Graphs of Polyhedra}
\address Institute
of Mathematics, Hebrew University of Jerusalem, Jerusalem, 
Israel and IBM Almaden Research Center, San Jose, 
California 95120\endaddress
\address Department of Mathematics, Massachusetts 
Institute of Technology,
Cambridge, Massachusetts 02139\endaddress
\date July 1, 1991\enddate
\thanks The first author was supported in part by a BSF 
grant by a GIF grant.
The second author was supported by an AFOSR grant\endthanks
\abstract The diameter of the graph of a $d$-dimensional 
polyhedron with
$n$ facets is at most $n^{\log d+2}$\endabstract
\endtopmatter

\document

Let $P$ be a convex polyhedron. The {\it graph\/} of $P$ 
denoted by $\roman G(P)$ is an
abstract graph whose vertices are the extreme points of 
$P$ and two vertices
$u$ and $v$ are adjacent if the interval $[v,u]$ is an 
extreme edge
$(=1$-dimensional face) of $P$. The diameter of the graph 
of $P$ is denoted by
$\delta(P)$. 

Let $\Delta(d,n)$ be the maximal diameter of the graphs of 
$d$-dimen-\linebreak
sional polyhedra $P$ with $n$ facets. (A facet is a 
$(d-1)$-dimensional face.) Thus,
$P$ is the set of solutions of $n$ linear inequalities in 
$d$ variables. It is
an old standing problem to determine the behavior of the 
function
$\Delta(d,n)$. The value of $\Delta(d,n)$ is a lower bound 
for the number of
iterations needed for Dantzig's simplex algorithm for 
linear programming with
any pivot rule.

In 1957 Hirsch conjectured \cite2 that $\Delta(d,n)\leq 
n-d$. Klee and Walkup
\cite6 showed that the Hirsch conjecture is false for 
unbounded polyhedra. They
proved that for $n\ge 2d, \Delta(d,n)\ge n-d+[d/5]$. This 
is the best
known lower bound for $\Delta(d,n)$.
 The
statement of the Hirsch conjecture for bounded polyhedra 
is still open. For a
recent survey on the Hirsch conjecture and its relatives, 
see \cite5.

In 1967 Barnette proved \cite{1, 3} that $\Delta(d,n)\leq 
n3^{d-3}$. An
improved upper bound, $\Delta(d,n)\leq n2^{d-3}$, was 
proved in 1970 by Larman
\cite7. Barnette's and Larman's bounds are linear in $n$ 
but exponential in the
dimension $d$. In 1990 the first author \cite4 proved a 
subexponential bound
$\Delta(d,n)\leq2^{\sqrt{(n-d)\log(n-d)}}$.

The purpose of this paper is to announce and to give a 
complete proof of a
quasi-polynomial upper bound for $\Delta(d,n)$. Such a 
bound was proved by the
first author in March 1991. The proof presented here is a 
substantial
simplification that was subsequently found by the second 
author. See \cite4
for the original proof and related results. The existence 
of a polynomial (or
even linear) upper bound for $\Delta(d,n)$ is still open.
Recently, the first author found a randomized pivot rule 
for linear programming
which requires an expected $n^{4\sqrt d}$ (or less) 
arithmetic operations
for every linear programming problem with $d$ variables 
and $n$ constraints.

\thm{Theorem 1} 
$$
\Delta(d,n)\leq n^{\log\,d+2}.
\tag1
$$
\ethm
\demo{Proof} Let $P$ be a $d$-dimensional polyhedron with 
$n$ facets, and let
$v$ and $u$ be two vertices of $P$. Let $k_v$ $[k_u]$ be 
the maximal positive
number such that the union of all vertices in all paths in 
$\roman G(P)$
starting from $v\ [u]$ of length at most $k_v$ $[k_u]$ are 
incident to at most
$n/2$ facets. Clearly, there is a facet $F$ of $P$ so that 
we can reach $F$ by
a path of length $k_v+1$ from $v$ and a path of length 
$k_u+1$ from $u$. We
claim now that $k_v\leq\Delta(d,[n/2])$. Indeed, let $Q$ 
be the polyhedron
obtained from $P$ by ignoring all the inequalities that 
correspond to facets that
cannot be reached from $v$ by a path of length at most 
$k_v$. Let $\omega$ be a
vertex in $\roman G(P)$ whose distance from $v$ is $k_v$. 
We claim that the
distance of $\omega$ from $v$ in $G(Q)$ is also $k_v$. To 
see this consider the
shortest path between $v$ and $\omega$ in $G(Q)$. If the 
length of this path is
smaller than $k_v$ there must be an edge in the path that 
is not an edge of
$P$. Consider the first such edge $E$. Since $E$ is not an 
edge of $P$, it
intersects a hyperplane $H$ that corresponds to one of the 
inequalities that
was ignored. This gives a path in $P$ of length smaller 
than $k_v$ from $v$ to
the facet of $P$ determined by $H$, which is a 
contradiction. \enddemo

We obtained that the distance from $v$ to $u$ is at most
$\Delta(d-1,n-1)$$+2\Delta(d,[n/2])+2$. This gives the 
inequality
$\Delta(d,n)\leq\Delta(d-1,n-1)$ $+ 2\Delta(d,[n/2])+2$, 
which implies the
statement of the theorem.

\enddocument